\documentclass[12pt]{amsart}
\usepackage{amssymb,amsmath,amsthm,amscd}
\usepackage[all]{xy}

\usepackage[hypertex]{hyperref}

\numberwithin{equation}{section}

\newtheoremstyle{fancy1}{10pt}{10pt}{\itshape}{12pt}{\textsc\bgroup}{.\egroup}{8pt}{
}
\newtheoremstyle{fancy2}{10pt}{10pt}{}{12pt}{\itshape}{.}{8pt}{ }

\theoremstyle{fancy1}

\newtheorem*{lem*}{Lemma}

\newtheorem{thm}[equation]{Theorem}
\newtheorem*{thm*}{Theorem}

\newtheorem*{main*}{Theorem}

\newtheorem*{cor*}{Corollary}
\newtheorem*{prop*}{Proposition}

\newtheorem*{problem*}{Problem}

\theoremstyle{fancy2}

\newtheorem*{rems*}{Remarks}

\newtheorem*{rem*}{Remark}

\newtheorem*{example*}{Example}

\newcommand{\cref}[1]{Corollary~\ref{#1}}

\newcommand{\e}{\epsilon}
\newcommand{\gs}{\sigma}

\newcommand{\RP}{\mathbb{R\mkern1mu P}}
\newcommand{\CP}{\mathbb{C\mkern1mu P}}
\newcommand{\HP}{\mathbb{H\mkern1mu P}}
\newcommand{\CaP}{\mathrm{Ca}\mathbb{\mkern1mu P}^2}
\newcommand{\Sph}{\mathbb{S}}
\newcommand{\Disc}{\mathbb{D}}

\newcommand{\C}{{\mathbb{C}}}
\newcommand{\R}{{\mathbb{R}}}
\newcommand{\Z}{{\mathbb{Z}}}

\renewcommand{\H}{\ensuremath{\operatorname{H}}}

\newcommand{\F}{\ensuremath{\operatorname{F}}}
\newcommand{\G}{\ensuremath{\operatorname{G}}}

\newcommand{\SO}{\ensuremath{\operatorname{SO}}}

\renewcommand{\O}{\ensuremath{\operatorname{O}}}
\newcommand{\Sp}{\ensuremath{\operatorname{Sp}}}
\newcommand{\U}{\ensuremath{\operatorname{U}}}
\newcommand{\SU}{\ensuremath{\operatorname{SU}}}
\newcommand{\Spin}{\ensuremath{\operatorname{Spin}}}

\newcommand{\T}{\ensuremath{\operatorname{T}}}
\renewcommand{\S}{\ensuremath{\operatorname{S}}}

\newcommand{\K}{\ensuremath{\operatorname{K}}}

\newcommand{\fg}{{\mathfrak{g}}}

\newcommand{\ft}{{\mathfrak{t}}}

\def\con#1=#2(#3){#1 \equiv #2 \bmod{#3}}

\newcommand{\ml}{\langle}                     
\newcommand{\mr}{\rangle}                    

\newcommand{\diag}{\ensuremath{\operatorname{diag}}}

\renewcommand{\sec}{\ensuremath{\operatorname{sec}}}
\newcommand{\Ric}{\ensuremath{\operatorname{Ric}}}

\DeclareMathOperator{\Isom}{Isom}

\newcommand{\Kpm}{K^{\scriptscriptstyle{\pm}}}
\newcommand{\Kp}{K^{\scriptscriptstyle{+}}}
\newcommand{\Km}{K^{\scriptscriptstyle{-}}}

\newcommand{\no}{\noindent}
\newcommand{\co}{{cohomogeneity}}
\newcommand{\coo}{{cohomogeneity one}}

\newcommand{\nnc}{non-negative curvature}
\newcommand{\nnsc}{non-negative sectional curvature}
\newcommand{\pc}{positive curvature}
\newcommand{\psc}{positive sectional curvature}

\newcommand{\hR}{\hat{R}}

\begin{document}

\title{ Riemannian Manifolds with positive sectional curvature}

\author{Wolfgang Ziller}
\address{University of Pennsylvania\\
   Philadelphia, PA 19104}
\email{wziller@math.upenn.edu}

\thanks{These are notes from a series of lectures given in Guanajuato, Mexico in 2010. The  author was supported by  a grant from the
National Science Foundation and by the Mexican National Academy of
Sciences.}

 \maketitle

 Of special interest in the history of Riemannian geometry have been
 manifolds with positive sectional curvature. In these notes we want
 to give a survey of this subject and some recent developments. We
 start with some historical developments.

\section{History and Obstructions}

It is fair to say that Riemannian geometry started with
 Gauss's famous "Disquisitiones generales" from 1827 in which one finds a
 rigorous discussion of what we now call the Gauss curvature of a
 surface. Much has been written
about the importance and influence of this paper, see in particular
the article \cite{Do} by P.Dombrowski for a careful discussion of
its contents and influence during that time. Here we only make a few
comments. Curvature of surfaces in 3-space had been studied
previously by a number of authors and was defined as the product of
the principal curvatures. But Gauss was the first to make the
surprising discovery that this curvature only depends on the
intrinsic metric and not on the embedding. Here one finds for
example the formula for the metric in
 the form $ds^2=dr^2+f(r,\theta)^2d\theta^2$. Gauss showed that every
 metric on a surface has this form in "normal" coordinates and
 that it has curvature $K=-f_{rr}/f$. In fact one can take it as the definition of the Gauss
 curvature and proves Gauss's famous "Theorema Egregium" that the curvature is an intrinsic invariant and does not depend on the embedding in $\R^3$.
  He also proved a
 local version of what we nowadays call the Gauss-Bonnet theorem (it
 is not clear what Bonnet's contribution was to this result), which
 states that in a geodesic triangle $\Delta$ with angles
 $\alpha,\beta,\gamma$ the Gauss curvature  measures the
 angle "defect":
 $$
\int_\Delta K dvol = \alpha+\beta+\gamma-\pi
 $$

 Nowadays the Gauss Bonnet theorem also goes under its global
 formulation for a compact surface:
$$
\int_M K dvol = 2\pi \chi(M)
 $$
 where $\chi(M)$ is the Euler characteristic. This follows from the defect formula by using a triangulation, but it is actually not
 found in any of Gauss's papers. Of course no rigorous definition of a manifold or of the Euler characteristic existed at the time.
  Maybe the first time the above  formulation can be
 found is in Blaschke's famous book "Vorlesungen ueber Differential
 Geometrie" from 1921 \cite{Bl} (although it is already discussed in a paper
 by Boy in 1903 \cite{Bo}).

 In any case, the formula implies that a compact surface with
 positive curvature must be the 2-sphere, or the real projective
 plane. This is of course the beginning of the topic of these
 lectures on  positive curvature.

The next big step was made by Riemann in his famous Habilitation
from 1854, eight months before Gauss's death). He started what we
now aptly call Riemannian geometry (in dimension bigger than $2$) by
giving intrinsic definitions of what is now called sectional
curvature (we will use $\sec$ for this notion instead of the more
common one $K$). For each 2-plane $\sigma\subset T_pM$ one
associates the sectional curvature, which we denote by
$\sec(\sigma)$. This can be defined for example as the Gauss
curvature of the 2-dimensional surface spanned by going along
geodesics in the direction of $\sigma$ (this was in fact one of
Riemann's definitions). Here one also finds for the first time an
explicit formula for a space of constant curvature $c$:
$$
ds^2=\frac{dx_1^2+\dots dx_n^2}{1+\frac{c}{4}(x_1^2+\dots x_n^2)}
$$
including in particular the important case of the hyperbolic plane
$c=-1$.

For our story, the next important development was Clifford's
discovery in 1873 of the Clifford torus
$\Sph^1(1)\times\Sph^1(1)\subset \Sph^3(\sqrt{2})\subset \R^4$,
which to his surprise has intrinsic  curvature $0$ (After all
something that looks like a plane has to extend to infinity). This
motivated Klein to formulate his famous Clifford-Klein space form
problem, which in one formulation asks to classify surfaces of
constant curvature. This has a painful history (after all, one needs
a good definition of completeness, a concept of a global surface and
some understanding of covering space theory). In a beautiful paper
by H.Hopf from 1926 \cite{Ho} he gave us our present definition of
completeness and solves the classification problem. It is amusing to
note that Hopf  points out that the many previous papers on the
subject, especially by Killing in 1891-1893, the authors did not
realize that the Moebius band has a flat metric.

For us the next development is of course the Bonnet-Myers theorem,
which holds more generally for the positivity of an average of the
sectional curvatures \cite{My2}:
\begin{thm*}[Bonnet-Myers] If $M$ has a complete metric with
$\Ric\ge 1$ then the diameter is at most $\pi$, and the fundamental
group is finite.
\end{thm*}
This theorem also has an interesting history. Bonnet in 1855 only
showed that the "extrinsic" diameter in 3-space has length at most
$\pi$. The difficulty to obtain an intrinsic proof in higher
dimensions was partially due to the fact that one needs a good
formula for the second variation, which surprisingly took a long
time to develop. Noteworthy are papers by Synge from 1925 \cite{Sy1}
(who was the first one to show that a geodesic of length $>\pi$
cannot be shortest by a second variation argument), Hopf-Rinow from
1931 \cite{HR} (where they proved any two points can be joined by a
minimal geodesic) Schoenberg from 1932 \cite{Sch}, Myers from 1935
\cite{My1} (here one finds for the first time the conclusion that
$\pi_1(M)$ is finite) and Synge \cite{Sy2} from 1935 as well. There
was a fierce competition between Myers and Synge for priorities
(\cite{My1} and \cite{Sy2} appeared in the same issue of Duke Math
J. and in Myer's paper one finds the mysterious footnote ``Received
by the Editors of the Annals of Mathematics, February 27, 1934,
accepted by them, and later transferred to this journal").
Schoenberg's paper contains the formula for second variation that
one now finds in books, and Synge's papers the usual proof in the
case of sectional curvature. In 1941 Myers used this proof and
summed over an orthonormal basis. Thus it would be fairer to call it
the Bonnet-Synge-Myers theorem. Nevertheless, Myers paper created a
lot of excitement at the time due to the importance of Ricci
curvature in general relativity.

Important for our story is another paper by Synge from 1936
\cite{Sy3} where he proved:
\begin{thm*}[Synge] If $M$ is a compact manifold with positive
sectional curvature, then $\pi_1(M)$ is $0$ or $\Z_2$ if $n$ is
even, and $M$ is orientable if $n$ is odd.
\end{thm*}
In particular, $\RP^n\times\RP^n$ does not admit a metric with
positive curvature. I can also recommend reading Preissman's paper
from 1936 \cite{Pr} on negative curvature, still very readable for
today's audience.

The surprising fact is that the above two theorems are the only
known obstructions that deal with positive curvature only. There are
a number theorems that give obstructions to non-negative curvature.
On the other hand, one expects that the class of manifolds admitting
positive curvature is much smaller than the class admitting
non-negative curvature (and this is born out in known examples).
Since this is not the purpose of the present notes, we just
summarize them:

\begin{itemize}
\item
(Gromov) If $M^n$ is a compact manifold with $\sec\ge 0$, then there
exists a universal constant $c(n)$ such that $b_i(M^n,F)\le c(n)$
for all $i$ and any field of coefficients $F$. Furthermore, the
fundamental group has a generating set with at most $c(n)$ elements.

\item (Cheeger-Gromoll)
If $M^n$ is a compact manifold that admits a metric with
non-negative sectional curvature,  then there exists an abelian
subgroup of $\pi_1(M^n)$ with finite index.

\item (Lichnerowicz-Hitchin) The obstructions to positive scalar curvature
imply that a compact spin manifold with $\hat{A}(M)\ne 0$ or
$\alpha(M)\ne 0$ does not admit a metric with non-negative
sectional curvature, unless it is flat. In particular there
exist exotic spheres, e.g. in dimension 9, which do not admit
positive curvature.

\item
(Cheeger-Gromoll) If $M^n$ is a non-compact manifold with a complete
metric with $\sec\ge 0$, then there exists a totally geodesic
compact submanifold $S^k$, called the soul, such that $M^n$ is
diffeomorphic to the normal bundle of $S^k$.
\end{itemize}

If we allow ourselves to add an upper as well as a lower bound on
the sectional curvature it is convenient to introduce what is called
the {\it pinching constant} which is defined as $\delta = \min \sec
/ \max\sec$. One then has the following recognition and finiteness
theorems:

\begin{itemize}
\item
(Berger-Klingenberg, Brendle-Schoen) If $M^n$ is a compact
manifold with $\delta\ge \frac 1 4 $, then $M$ is either
diffeomorphic to a space form $\Sph^n/\Gamma$ or isometric to
$\CP^n$, $\HP^n$ or $\CaP$ with their standard Fubini metric.
\item
(Cheeger) Given a positive constant $\e$, there are only finitely
many diffeomorphism types of compact simply connected manifolds
$M^{2n}$ with $\delta\ge \e$.
\item (Fang-Rong,Petrunin-Tuschmann) Given a positive constant $\e$, there are only finitely many
diffeomorphism types of compact   manifolds $M^{2n+1}$ with
$\pi_1(M)=\pi_2(M)=0$ and $\delta \ge \e$.
\end{itemize}

Since our emphasis is positive curvature, we will not discuss other
results, except in passing, about non-negative curvature in this
survey. We finally mention some conjectures.

\begin{itemize}
\item
(Hopf) There exists no metric with positive sectional curvature on
$\Sph^2\times\Sph^2$. More generally, there are no positively curved
metrics on the product of two compact manifolds, or on a symmetric
space of rank at least two.
\item (Hopf) A compact manifold with $\sec\ge 0$ has non-negative Euler
characteristic. An even dimensional manifold with positive curvature
has positive Euler characteristic.
\item (Bott-Grove-Halperin) A
compact simply connected manifold $M$ with $\sec\ge 0$ is {\it
elliptic}, i.e.,  the sequence of Betti numbers of the loop space of
$M$
 grows at
 most polynomially for every field of coefficients.
\end{itemize}
The latter conjecture , and its many consequences, were discussed in
the literature for the first time in \cite{GH}. It is usually
formulated for rational coefficients, where it is equivalent to the
condition that only finitely many homotopy groups are non-zero
(called rationally elliptic). One can thus apply rational homotopy
theory to obtain many consequences. E.g., it implies, under the
assumption of non-negative curvature, that $b_i(M^n,F)\le 2^n$ and
that the Euler characteristic is non-negative (Hopf conjecture), and
positive in even dimensions iff all odd Betti numbers vanish. The
above more geometric formulation, which one should call elliptic, is
a natural generalization. If $n=4$, rational homotopy theory implies
that $M$, if compact and simply connected, is diffeomorphic to one
of the known examples of non-negative curvature, i.e.
$\Sph^4,\CP^2,\Sph^2\times\Sph^2$ or $\CP^2\# \pm \CP^2$.  In
\cite{PP} it was shown that a compact  simply connected elliptic
5-manifold is diffeomorphic to one of the known examples with \nnc,
i.e., one of $\Sph^5$, $\SU(3)/\SO(3)$, $\Sph^3\times\Sph^2$ or the
non-trivial $\Sph^3$ bundle over $\Sph^2$. In both cases, no
curvature assumption is necessary.

\smallskip

 Of course,
one should  also  mentions Hamilton's theorem which states that a 3
manifold with positive curvature is diffeomorphic to a space form
$\Sph^4/\Gamma$. Thus in dimension 2 and 3, manifolds with positive
curvature are classified.

\bigskip

We formulate some other natural conjecturs:
\begin{itemize}
\item A compact simply connected 4 manifold with \pc\ is
diffeomorphic to $\Sph^4$ or $\CP^2$.
\item  A compact simply connected 5 manifold with \pc\ is
diffeomorphic to $\Sph^5$.
\item (Klingenberg-Sakai) There are only finitely many diffeomorphism classes of
positively curved manifolds in a given homotopy type.
\item There are only finitely many diffeomorphism classes of
positively curved manifolds in even dimensions, and all odd
Betti numbers are 0.
\item In odd dimension, there are only finitely many 2-connected
manifolds with positive curvature.
\end{itemize}

The last 2 finiteness conjectures are probably too optimistic, but
one should at least expect an upper bound on the Betti numbers, e.g.
at most 2 in dimension 6.

\section{Compact examples of positive curvature}

Homogeneous spaces which admit a homogeneous metric with \pc\ have
been classified by Wallach in even dimensions (\cite{Wa}) and by
B\'erard-Bergery in odd dimensions (\cite{BB}). We now describe
these examples, due to Berger, Wallach and Aloff-Wallach
\cite{Be,Wa,AW}, as well as the biquotient examples \cite{E1, E2,
Ba}. In most cases we will also mention that they admit natural
fibrations, a topic we will cover in Section 4.
\smallskip

1) The god given basic examples of positive curvature are the rank
one symmetric spaces $\Sph^n$, $\CP^n$, $\HP^n$ or $\CaP$. (We do
not know where in the literature  it is first discussed that $\CaP$
carries a metric with positive curvature which is $1/4$ pinched).
They admit the well known homogeneous Hopf fibrations. Recall that a
homogeneous  fibration is of the form $\K/\H\to\G/\H\to\G/\K$
obtained from inclusions $\H\subset\K\subset\G$.

\bigskip

 \centerline{ $\Sph^1\to \Sph^{2n+1}\to \CP^n$ obtained from  $\SU(n)\subset
    \U(n)\subset \SU(n+1)$,}

\bigskip

 \centerline{ $\Sph^3\to \Sph^{4n+3}\to \HP^n$  obtained from  $\Sp(n)\subset
    \Sp(n)Sp(1)\subset \Sp(n+1)$,}

\bigskip

 \centerline{ $\Sph^2\to\CP^{2n+1}\to\HP^n$  obtained from  $\Sp(n)U(1)\subset
    \Sp(n)\Sp(1)\subset \Sp(n+1)$.}

\bigskip

 \centerline{ $\Sph^7\to \Sph^{15}\to \Sph^8$ coming from $\Spin(7)\subset
    \Spin(8)\subset \Spin(9)$.}

\bigskip

2) The homogeneous flag manifolds due to Wallach: $W^6=\SU(3)/\T^2$,
$W^{12}=\Sp(3)/\Sp(1)^3$ and $W^{24}=\F_4/\Spin(8)$. They are the
total space of the following homogeneous fibrations:
\smallskip

 \centerline{$\Sph^2\to \SU(3)/\T^2 \to \CP^2$,}

\bigskip

 \centerline{$\Sph^4\to \Sp(3)/\Sp(1)^3 \to \HP^2$,}

\bigskip

 \centerline{$\Sph^8\to \F_4/\Spin(8) \to \CaP$.}

\bigskip

3) The Berger space $B^{13}=\SU(5)/\Sp(2)\cdot \S^1$  admits a
fibration

\bigskip
 \centerline{$\RP^5\to\SU(5)/\Sp(2)\cdot \S^1\to\CP^4,$}
\bigskip

\no coming from the inclusions $\Sp(2)\cdot
\S^1\subset\U(4)\subset\SU(5)$. Here $\Sp(2)\subset\SU(4)$ is the
usual embedding and $\S^1$ is the center of $\U(4)$. Furthermore,
the fiber is $\U(4)/\Sp(2)\cdot \S^1=\SU(4)/\Sp(2)\cdot \Z_2
=\SO(6)/\O(5)=\RP^5$.
\bigskip

4) The Aloff-Wallach spaces
$W^7_{p,q}=\SU(3)/\diag(z^{p},z^{q},\bar{z}^{p+q})\ , \gcd(p,q)=1$.
By interchanging coordinates we can assume $p\ge q\ge 0$. They have
positive curvature, unless $(p,q)=(1,0)$. They also admit
interesting fibrations

\bigskip
\centerline{$\Sph^3/\Z_{p+q}\to W_{p,q}\to \SU(3)/\U(2),$}
\bigskip

\no coming from the inclusions
$\diag(z^{p},z^{q},\bar{z}^{p+q})\subset\U(2)\subset\SU(3)$. Hence,
as long as $p+q\ne 0$ (or $q\ne 0$ in the above notation), the fiber
is the lens space
$\U(2)/\diag(z^{p},z^{q})=\SU(2)/\diag(z^{p},z^{q})$ with
$z^{p+q}=1$. In the special case of $p=q=1$, we obtain a principal
$\SO(3)$ bundle.

\smallskip

Another fibration is of the form

\bigskip
\centerline{$\Sph^1\to W_{p,q}\to \SU(3)/\T^2,$}
\bigskip
\no coming from the inclusions
$\diag(z^{p},z^{q},\bar{z}^{p+q})\subset\T^2\subset\SU(3)$.

\bigskip

5) The Berger space: $B^7=\SO(5)/\SO(3)$. To describe the embedding
$\SO(3)\subset\SO(5)$, we recall that $\SO(3)$ acts orthogonally via
conjugation on the set of $3\times 3$ symmetric traceless matrices.
This space is special since $\SO(3)$ is maximal in $\SO(5)$ and
hence does not admit a homogeneous fibration. On the other hand, in
\cite{GKS} it was shown that the manifold is diffeomorphic to an
$\Sph^3$ bundle over $\Sph^4$. It is also what is called isotropy
irreducible, i.e., the isotropy action of $H$ on the tangent space
is irreducible. This implies that there is only one $\SO(5)$
invariant metric up to scaling.

\bigskip

Thus  all of these examples in 2)- 5) are  the total space of a
fibration. This property will be interesting to us in Section 4.

\bigskip

A natural generalization of homogeneous spaces are so called
biquotients, discussed for the first time in \cite{GM}. For this,
let $\G/\H$ be a homogeneous space and $\K\subset \G$ a subgroup.
Then $\K$ acts on $\G/\H$ on the left, and in some cases the action
is free, in which case the manifold $\K\backslash \G/\H$ is a
biquotient. An equivalent formulation is as follows: Take a subgroup
$\U\subset \G\times\G$ and let $\U$ act one the left and right
$(u_1,u_2)* g=u_1gu_2^{-1}$. The action is free, if for any
$(u_1,u_2)\in \U$ the element $u_1$ is not conjugate to $u_2$ unless
$u_1=u_2$ lies in the center of $\G$. We denote the quotient by
$\G/\!/\U$. The biinvariant metric on $G$ (or $\G\times\G$ ) induces
a metric on $G/\!/U$ with \nnsc. In some cases, this can be deformed
(via a Cheeger deformation) into one with \pc.   We now describe
these biquotient examples, due to Eschenburg and Bazaikin,
explicitly.

\smallskip

6) There is an analogue of the 6-dimensional flag manifold which is
a biquotient of $\SU(3)$ under an action of $\T^2=\{(z,w)\mid
z,w\in\C \ , \ |z|=|w|=1\}$. It is given by:
$$E^6=\SU(3)/\!/\T^2=\diag(z,w,zw)\backslash\SU(3)/\diag(1,1,z^2w^2)^{-1}.$$
The action by $\T^2$ is clearly free. In order to show that  this
manifold is not diffeomorphic  to the homogeneous flag $W^6$, one
needs to compute the cohomology with integer coefficients. The
cohomology groups are the same for both manifolds, but the ring
structure is different (\cite{E2}). The examples $W^6$ and $E^6$,
which have $b_2=2$, as well as $\Sph^6$ and $\CP^3$, are the only
known examples of \pc. It is thus a natural question wether \pc\ in
dimension 6 implies that the Betti numbers satisfy $b_1=b_3=b_5=0$
and $b_2=b_4\le 2$.

The inhomogeneous flag also admits a fibration of a (different)
sphere bundle similar to the flag manifold:

\smallskip
 \centerline{$\Sph^2 \to \SU(3)/\!/\T^2 \to \CP^2$}
\smallskip

\bigskip

7) We now describe the 7-dimensional family of Eschenburg spaces
$E_{k,l}$, which can be considered as a generalization of the Aloff
Wallach spaces. Let $k:= (k_1,k_2,k_3)$ and $l:= (l_1,l_2,l_3) \in
\Z^3$ be two triples of integers with $\sum k_i=\sum l_i $. We can
then define a two-sided action of $\S^1 = \{z\in \C \mid |z|=1 \}$
on $\SU(3)$ whose quotient we denote by $E_{k,l}$:
$$E_{k,l} = \SU(3)/\!/\S^1=  \diag(z^{k_1}, z^{k_2},
z^{k_3})\backslash \SU(3)/ \diag(z^{l_1}, z^{l_2}, z^{l_3})^{-1}.$$
\no The action is free if and only if
$\diag(z^{k_1},z^{k_2},z^{k_3})$ is not conjugate to
$\diag(z^{l_1},z^{l_2},z^{l_3})$, i.e.
\begin{equation*}\label{free}
\gcd(k_1-l_{i}\ , k_{2}-l_{j})=1,\ \ \text{\it for all }\ \  i\neq
j\,,\,i,j \in \{1,2,3\}\,.
\end{equation*}

\no Eschenburg showed that $E_{k,l}$ has \pc\ if
$$
k_i \notin [\min(l_1,l_2,l_3),\max(l_1,l_2,l_3)].
$$

Among the biquotients $E_{k,l}$ there are two interesting
subfamilies. $E_p=E_{k,l}$ with $k=(1,1,p)$ and $l=(1,1,p+2)$ has
positive curvature when $p>0$. It admits a large group acting by
isometries. Indeed, $G=\SU(2)\times\SU(2)$ acting on $\SU(3)$ on the
left and on the right, acts by isometries in the Eschenburg metric
and commutes with the $\S^1$ action. Thus it acts by isometries on
$E_p$ and one easily sees that $E_p/G$ is one dimensional, i.e.,
$E_p$ is \coo. A second family consists of the \co\ two Eschenburg
spaces $E_{a,b,c}=E_{k,l}$ with $k=(a,b,c)$ and $l=(1,1,a+b+c)$.
Here $c=-(a+b)$ is the subfamily of Aloff-Wallach spaces. The action
is free iff $a,b,c$ are pairwise relatively prime and the Eschenburg
metric has \pc\ iff, up to permutations, $a\ge b\ge c >0$ or $a\ge b
>0\ , c<-a$.
For these spaces $G=\U(2)$ acts by isometries on the right and
$E_{a,b,c}/G$ is two dimensional. For a general Eschenburg space
$G=\T^3$ acts by isometries and $E_{k,l}/G$ is four dimensional. In
\cite{GSZ} it was shown that these groups $G$ are indeed the id
component of the full isometry group of a positively curved
Eschenburg space, unless it is an Aloff-Wallach space.

\smallskip
There are again natural fibrations. In the case of $E_{a,b,a+b}$
with $a\ge b>0$, the circle fibrations:

\smallskip
\centerline{$\Sph^1\to E_{a,b,a+b}\to \SU(3)/\!/\T^2,$}
\smallskip

\no and  the lens space fibrations:

\smallskip
\centerline{$\Sph^3/\Z_{a+b}\to E_{a,b,a+b}\to \CP^2$}
\smallskip

\no which, in the case of $a=b=1$, gives a second $\SO(3)$ principal
bundle over $\CP^2$.

\smallskip

The cohomogeneity two Eschenburg spaces admit orbifold fibrations,
which will also be of interest to us in Section 4.

\bigskip
\centerline{$F\to E_{a,b,c}\to \CP^2[a+b,a+c,b+c]\ ,$}
\bigskip

\no where the fiber $F$ is $\RP^3$ if all $a,b,c$'s are odd, and
$F=\Sph^3$ otherwise. Here the base is a \mbox{2-dimensional}
weighted complex projective space. A general Eschenburg space is the
total space of an orbifold circle bundle \cite{FZ1}.

\bigskip

8) We finally have the 13-dimensional Bazaikin spaces $B_{q}$, which
can be considered as a generalization of the Berger space $B^{13}$.
Let $q=(q_0,\dots,q_5)$ be a 6-tuple of integers with $\sum q_i=0$
and define
\begin{equation*}\label{Ba2}
B_q = \diag (z^{q_1} , \dots , z^{q_5} ) \backslash \SU(5)/ \diag
(z^{{-q_0}},A)^{-1},
\end{equation*}
where $A\in \Sp(2)\subset \SU(4)\subset\SU(5)$. Here we follow  the
treatment in \cite{Zi1} of Bazaikin's work \cite{Ba} (see also
\cite{EKS}). First, one easily shows that the action of
$\Sp(2)\cdot\S^1$ is free if and only if
\begin{equation*}\label{e:gcd}
\text{\it all $q_i$'s are odd\ \ and\ \ }\gcd(q_{\gs (1)}+q_{\gs
(2)},q_{\gs (3)}+ q_{\gs (4)})=2,
\end{equation*}
for all permutations $\gs\in S_5.$ On $\SU(5)$ we choose an
Eschenburg metric by scaling the biinvariant metric on $\SU(5)$  in
the direction of $\U(4)\subset\SU(5)$. The right action of
$\Sp(2)\cdot\S^1$ is then  by isometries. Repeating the same
arguments as in the previous case, one shows that the induced metric
on $\SU(5)/\!/\Sp(2)\cdot\S^1$ satisfies
\begin{equation*}\label{pos}
\sec>0 \text{\it \quad if and only if\quad }  q_i+q_j>0 \ (\text{or}
<0)\; \text{\it for all }  i < j .
\end{equation*}
The special case of $q=(1,1,1,1,1)$ is the homogeneous Berger space.
One again has a one parameter subfamily that is \coo, given by
$B_p=B_{(1,1,1,1,2p-1)}$ since $\U(4)$ acting on the left induces an
isometric action on the quotient. It has \pc\ when $p\ge 1$.

There is another equivalent description of the Bazaikin spaces given
by
\begin{equation*}\label{Ba2}
B_q = \diag (z^{q_1} , \dots  , z^{q_6}) \backslash \SU(6)/ \Sp(3)
\end{equation*}
with $\sum q_i=0$.

\smallskip

For these manifolds one has natural fibrations obtained from both
descriptions, given by

\smallskip \centerline{$\Sph^1\to \SU(6)/\Sp(3)\to B_q,$}
\smallskip

\no and

\centerline{$\Sph^5\to \SU(6)/\Sp(3)\to\Sph^9$.}
\smallskip
\no But in this case $B_q$ is not the total space of a fibration,
unless it is homogeneous. On the other hand, if we allow orbifold
fibrations they all admit one:

\bigskip \centerline{$\RP^5\to B_q\to
\CP^4[q_0+q_1,q_0+q_2,\dots,q_0+q_6].$}
\smallskip

\bigskip

Unlike in the homogeneous case, there is no general classification
of positively curved biquotients, except  in the following cases. We
call a  metric on $G /\!/H$ torus invariant if it is induced by a
left invariant metric on $G$ which  is also right invariant under
the action of a maximal torus. The main theorem in \cite{E2} states
that an even dimensional biquotient $G /\!/H$ with $G$ simple and
which admits a positively curved torus invariant metric is
diffeomorphic to a rank one symmetric space or the biquotient
$\SU(3)/\!/T^2$. In the odd dimensional case he shows that $G /\!/H$
with a positively curved torus invariant metric and $G$ of rank 2 is
either diffeomorphic to a homogeneous space or a positively curved
Eschenburg space.

\bigskip

There are only two more examples which are not homogeneous or
biquotients. One is a 7-dimensional exotic sphere due to
Petersen-Wilhem \cite{PW} (although the rather delicate calculations
have not yet been verified). The method is via deforming a natural
metric of non-negative curvature on a biquotient description of the
exotic sphere to positive curvature. Deforming \nnc\ to \pc\ is an
important problem, and not yet well understood.

\smallskip

The second example is due to Grove-Verdiani-Ziller \cite{GVZ}, and
independently O.Dearricott \cite{De}, and will be discussed in
Section 4.
 It arises as the total space of an orbifold fibration.

 \bigskip It is
also interesting to examine the topology of the known examples. In
\cite{KS} it was shown that there exist pairs of positively curved
Aloff Wallach spaces which are homeomorphic but not diffeomorphic.
This turns out to happen more frequently for the Eschenburg spaces
\cite{CEZ}. For such pairs $M, M'$ one knows that $M =
M'\#\Sigma^7$, for some exotic sphere $\Sigma^7$. It is not hard to
check that among the examples in \cite{CEZ}, every exotic 7-sphere
can occur as a factor $\Sigma^7$, whereas this does not seem to be
the case for Aloff Wallach spaces. On the other hand, \cite{FZ2}
provides evidence that positively curved Bazaikin spaces are
homeomorphically distinct.

\bigskip

\section{Positive curvature with symmetry}

As we saw in Section 1, not much is known as far as general
obstructions to positive curvature is concerned. A very successful
program was suggested by K.Grove, motivated by the Kleiner-Hsiang
theorem below \cite{HK}, that one should examine positive curvature
under the additional assumption of a large symmetry group.
\begin{thm*}[Kleiner-Hsiang] If $M$ is a compact simply connected
4-manifold on which a circle acts by isometries, then $M$ is
homeomorphic to $\Sph^4$ or $\Sph^2$.
\end{thm*}
Topological results on circle actions imply that they are
diffeomorphic, and a recent result by Grove-Wilking shows that the
$\S^1$ action is linear.
 Thus a counter example to the Hopf conjecture would have to have a
 finite isometry group.

 In higher dimensions, one obtains obstructions assuming that a
 torus of large dimension acts \cite{GS}.
\begin{thm*}[Grove-Searle] If $M^n$ is a compact simply connected
manifold with positive curvature on which a torus $\T^s$ acts by
isometries, then $s\le n/2$ in even dimensions, and $s\le (n+1)/2$
in odd dimensions. Equality holds iff M is diffeomorphic to $\Sph^n$
or $\CP^n$.
\end{thm*}
See the article in this volume where it is shown that the torus
action is linear as well. Great strides were made by B.Wilking
\cite{Wi1,Wi2} who showed
\begin{thm*}[Wilking] Let $M^n$ be a compact simply connected
manifold with \pc.
\begin{enumerate}
\item[(a)] If $n\ne 7$ and $\T^s$ acts by isometries with $s\ge\frac{n+1}{4}$, then
$M$ is homotopy equivalent to a rank one symmetric space.
  \item[(b)] If $0<\dim M/G<  \sqrt{n/18}-1$, then $M$ is homotopy equivalent to a rank one symmetric space.
   \item[(c)] If $\dim \Isom(M)\ge 2n-6 $ then $M$ is either homotopy equivalent to a rank one symmetric
   space or isometric to a homogeneous space with positive
   curvature.
\end{enumerate}
\end{thm*}

One of the main new tools is the so called connectedness Lemma,
which turns out to be very powerful.
\begin{lem*}[Connectedness Lemma]
If $M^n$ has positive curvature, and $N$ is a totally geodesic
submanifold of codimension $k\le (n+1)/2$, then the inclusion $N
\hookrightarrow M$ is $n-2k+1$ connected.
\end{lem*}
The proof is surprisingly simple and similar to the proof of Synge's
theorem. One shows that in the loop space $\Omega_N(M)$ of curves
starting and ending at $N$ every critical point, i.e. geodesic
starting and ending perpendicular to $N$, has index at least
$n-2k+1$ since there are $n-2k+1$ parallel Jacobi fields starting
and ending perpendicular to $N$. This implies that the inclusion
$N\to \Omega_N(M) $ is $n-2k$ connected, and hence $N\to M$ is
$n-2k+1$ connected.

 An important consequence is a certain kind of periodicity
in cohomology:
\begin{lem*}[Periodicity Theorem]
If $M^n$ has positive curvature, and $N$ is a totally geodesic
submanifold of codimension $k$, then there exists cohomology class
$e\in H^k(M,\Z)$ such that $\cup e :H^i(M,\Z)\to H^{i+k}(M,\Z)$ is
an isomorphism for $k\le i \le n-2k$.
\end{lem*}

Along the way to proving these results, he obtains a number of
fundamental obstructions  on the structure of the possible isotropy
groups of the action. We mention a few striking examples.
\begin{thm*}[Wilking] Let $M^n$ be a compact simply connected
manifold with \pc\ on which $G$ acts by isometries with principal
isotropy group $H$.
\begin{enumerate}
\item[(a)] If $H$ is non-trivial, then $\partial(M/G)$ is non-empty.
  \item[(b)] Every irreducible subrepresentation of $G/H$ is a
  subrepresentation of  $K/H$  where $K$ is an isotropy group
  and $K/H$ is a sphere.
   \item[(c)] If $\dim (M/G)=k$, then $\partial(M/G)$ has at most
   $k+1$ faces, and in the case of equality $M/G$ is
   homeomorphic to a simplex.
\end{enumerate}
\end{thm*}
Part (a) is powerful since the distance to the boundary is a
strictly convex function.  In general one can use Alexandrov
geometry on the quotient as is an important tool, see the article by
Fernando Galaz-Garcia in this volume. Part (b) has strong
implications for  the pair $(G,H)$, with a very short list of
possibilities when the rank of $ H$ is bigger than 1.

\bigskip

Recently, L.Kennard proved two theorems concerning the Hopf
conjectures with symmetry \cite{Ke}, see also his article in this
volume.

\begin{thm*}[Kennard] Let $M^n$ be a compact simply connected
manifold with \psc.
\begin{enumerate}
\item[(a)] If $n=4k$ and $\T^r$ acts effectively and isometrically with $ r\ge
2\log_2(n)$, then $\chi(M) > 0$.
  \item[(b)] Suppose $M^n$ has the rational cohomology of a simply connected, com-
pact symmetric space $N$. If $\T^r$ acts  isometrically with $r
\ge 2 \log_2 n+7$, then $N$ is a product of spheres times either
a rank one symmetric space or a rank p Grassmannian
$\SO(p+q)=\SO(p)\SO(q)$ with $p=2$ or $p=3$.
\end{enumerate}
\end{thm*}

The main new tool is to use the action of the Steenrod algebra to
improve periodicity theorems. For example, the analogue of the
connectedness Lemma is
\begin{thm*}[Kennard] Let $M^n$ be a compact simply connected
manifold.
\begin{enumerate}
\item[(a)] If $M^n$ has \pc\ and contains a pair of totally geodesic, transversely
intersecting submanifolds of codimensions $k_1,k_2$ such that
$2k_1 + 2k_2 \le  n$, then $H^*(M;Q)$ is $\gcd(4, k_1,
k_2)$-periodic.
  \item[(b)] If $H^*(M;\Z)$ is
k-periodic with $3k\le n$, then $H^*(M;Q)$ is $\gcd(4,
k)$-periodic.
\end{enumerate}
\end{thm*}
Here the cohomology is called $k$-periodic if there exists
cohomology class $e\in H^k(M,\Z)$ such that $\cup e :H^i(M,\Z)\to
H^{i+k}(M,\Z)$ is an isomorphism for $0 < i < n - k$,  surjective
for $i=0$ and injective for $i=n-k$. In particular,
$H^{ik}(M,\Z)\simeq\Z$ for $0\le i\le n-2k-1 $.

\bigskip

One conclusion one can draw from these results is that positive
curvature with a large isometry group can only be expected in low
dimensions. This is born out in the classification of cohomogeneity
one manifolds with positive curvature. We first need to describe the
structure of such manifolds.

A simply connected compact cohomogeneity one manifold is  the union
of two homogeneous disc bundles. Given  compact Lie groups $H,\, \Km
,\, \Kp$ and $\G$ with inclusions $H\subset \Kpm \subset G$
satisfying $\Kpm/H=\Sph^{\ell_\pm}$, the transitive action of $\Kpm$
on
  $\Sph^{\ell_\pm}$ extends to a linear action on the disc $\Disc^{{\ell_\pm}+1} $.
We can thus define $$M=G\times_{\Km}\Disc^{{\ell_-}+1}\cup
G\times_{\Kp}\Disc^{{\ell_+}+1}$$ glued along the boundary $
\partial
(G\times_{\Kpm}\Disc^{\ell_\pm+1})=G\times_{\Kpm}\Kpm/H=G/H$ via the
identity. $G$  acts on $M$ on each half via left action in the first
component. This action has principal isotropy group $H$ and singular
isotropy groups $\Kpm$.
 One possible description of a cohomogeneity one manifold is thus
simply in terms of the Lie groups $H\subset \{\Km , \Kp\}\subset G$.

The simplest example is $\{e\}\subset \{\S^1 , \S^1\}\subset \S^1$
which is the manifold $\Sph^2$ with $G=\S^1$ fixing north and south
pole (and thus $\Kpm=G$) and principal isotropy trivial. The
isotropy groups $\{e\}\subset \{\S^1\times\{e\} ,\{e\}\times
\S^1\}\subset \S^1\times\S^1$ describe the 3-sphere $\Sph^3\subset
\C\oplus\C$ on which $\S^1\times\S^1 $ acts in each coordinate. More
subtle is the example $\S(\O(1)\O(1)\O(1))\subset \{\S(\O(2)\O(1)) ,
\S(\O(1)\O(2)\}\subset \SO(3)$. This is the 4-sphere, thought of as
the unit sphere in the set of $3\times 3$ symmetric traceless
matrices, on which $\SO(3)$ acts via conjugation.

\smallskip

The first new family of cohomogeneity one manifolds we denote by
$P_{ (p_-,q_-),(p_+,q_+)}$,
 and is given by the group diagram
$$
 H=\{\pm (1,1),\pm
(i,i),\pm (j,j),\pm (k,k)\} \subset\{ (e^{ip_-t},e^{iq_-t})\cdot H\;
   ,\;
   (e^{jp_+t},e^{jq_+t})\cdot H
   \}\subset\S^3\times\S^3.
  $$
  where $\gcd(p_-,q_-)=\gcd(p_+,q_+)=1$ and all $4$ integers are
   congruent to $1$ mod $4$.

  The second family  $Q_{
(p_-,q_-),(p_+,q_+)}$
  is given by the group diagram
$$
 H=  \{(\pm
1, \pm 1) , (\pm i , \pm i)\} \subset\{ (e^{ip_-t},e^{iq_-t})\cdot H\;
   ,\;
   (e^{jp_+t},e^{jq_+t})\cdot H
   \}\subset\S^3\times\S^3,
  $$
   where $\gcd(p_-,q_-)=\gcd(p_+,q_+)=1$, $q_+$ is even, and
   $p_-,q_-,p_+$ are
   congruent to $1$ mod $4$.

  \smallskip

  Special among these are the manifolds $P_k=P_{(1,1),(1+2k,1-2k)}$,
  $Q_k=Q_{(1,1),(k,k+1)}$ with $k\ge 1$, and the exceptional manifold $R^7=Q_{(-3,1),(1,2)}$. In terms of
these descriptions, we can state the classification, see [Ve, GWZ].

\begin{thm*}(Verdiani, $n$ even, Grove-Wilking-Ziller, $n$ odd)
A simply connected cohomogeneity one manifold $M^n$ with an
invariant metric of positive sectional curvature is equivariantly
diffeomorphic to one of the following:
\begin{itemize}
\item
An isometric action on a rank one symmetric space,
\item
     One of $E^7_p,  B^{13}_p$ or $B^7$,
\item
One of the  $7$-manifolds $P_k=P_{(1,1),(1+2k,1-2k)}$,
  $Q_k=Q_{(1,1),(k,k+1)}$ with $k\ge 1$, or the exceptional
  manifold $R^7=Q_{(-3,1),(1,2)}$
\end{itemize}
with one of the actions described above.
\end{thm*}
Here $P_k, Q_k,R$ should be considered as candidates for positive
curvature. Recently the exceptional manifold $R$ was excluded
\cite{VZ2}:

\begin{thm*}[Verdiani-Ziller]
Let $M$ be one of the $7$-manifolds $Q_{ (p_-,q_-),(p_+,q_+)}$ with
its cohomogeneity one action by $G=\S^3\times\S^3$ and assume that
$M$ is not of type  $Q_k,k\ge 0$. Then there exists no analytic
metric with non-negative sectional curvature invariant under $G$,
although there exists a smooth one.
  \end{thm*}
  In particular, it cannot carry an invariant metric with positive
  curvature.
  \smallskip

  Among the candidates $P_k,Q_k$ he first in each sequence  admit an invariant metric with
positive curvature since $P_1=\Sph^7$ and $Q_1=W^7_{1,1}$. The first
success in the Grove program to find a new example with positive
curvature is the following \cite{GVZ,De}:

\begin{thm*}[Grove-Verdiani-Ziller, Dearricott]
The cohomogeneity one manifold $P_2$ carries an invariant metric
with positive curvature.
  \end{thm*}
  As for the topology of this manifold one has the following
  classification \cite{Go}:

\begin{thm*}[Goette]
The cohomogeneity one manifold $P_k$ is diffeomorphic to
$E_k\#\Sigma^{\frac{k(k+1)}{2}}$, where $E_k$ is the $\S^3$
principal bundle over $\Sph^4$ with Euler class $k$, and $\Sigma$ is
the Gromoll-Meyer generator in the group of exotic 7-spheres.
  \end{thm*}
  In particular, it follows that $P_2$ is homeomorphic but not
  diffeomorphic to $T_1\Sph^4$, and thus indeed a new example of
  positive curvature.

\bigskip

\section{Fibrations with Positive Curvature}

  \smallskip

  As we saw in Section 2, many of the examples of positive curvature
  are the total space of a fibration. It is thus natural to ask
  under what condition the total space admits positive curvature,
  if  the base and fiber do. One should certainly expect conditions,
  since the bundle could be trivial.

  A.Weinstein examined this question in the context of Riemannian
  submersions with totally geodesic fibers \cite{We}. He called a bundle fat if $\sec(X,U)>0$ for all
  2-planes spanned by a vector $U$ tangent to the fibers and $X$ orthogonal to the fibers.
   For simplicity, we
  restrict ourselves for the moment to principle bundles. Let
  $\pi\colon
  P\to B$ be a $G$-principle bundle. Given a metric on the base $\ml . ,
  . \mr_B$, a principal connection $\theta\colon TP\to \fg$, and  a fixed
  biinvariant metric $Q$ on $G$, one defines
  a Kaluza Klein metric on $P$ as:
$$g_t(X,Y)=tQ(\theta(X),\theta(Y))+g(\pi_*(X),\pi_*(Y)). $$
Here one has the additional freedom to modify $t$, in fact $t\to 0$
usually increases the curvature.

The projection $\pi$ is then a Riemannian submersion with totally
geodesic fibers isometric to $(G,tQ)$.  Weinstein observed that the
fatness condition (for any $t$) is equivalent to requiring  that the
curvature $\Omega$ of $\theta$ has the property that
$\Omega_u=Q(\Omega,u)$ is a symplectic 2-form on the horizontal
space, i.e. the vector space orthogonal to the fibers, for every
$u\in\fg$. If $G=\S^1$, this is equivalent to the base being
symplectic. Fatness of the principal connection is already a strong
condition on the principal bundle which one can express in terms of
the characteristic classes of the bundle.

A.Weinstein made the following observation. Assume that $G$ and
$B^{2n}$ are compact and connected. For each  $y\in\fg$,   we have a
polynomial $q_y:\fg\to\R$ given by
\begin{equation*}\label{e:main}
q_y(\alpha)=\int_G\ml Ad_g(y),\alpha\mr^{2n} dg
\end{equation*}
which one easily checks is $Ad_G$--invariant (in $\alpha$ and $y$).
By Chern--Weyl theory, there exists a closed $2n$--form $\omega_y$
on $B^{2n}$ such that $\tau^*\omega_y=q_y(\Omega)$ and
$[\omega_y]\in H^{2n}(B,\R)$ represents a characteristic class of
the bundle. If the bundle is fat, $\Omega_y^{2n}\ne 0$ is a volume
form on $H$. Thus, if $G$ is connected, $\ml Ad_g(y),\Omega\mr^{2n}$
is nowhere zero and has constant sign when $g$ varies along $G$, and
the integral is thus $q_y(\Omega)$  nonzero on $H$. Hence $\omega_y$
is a volume form of $B^{2m}$, in particular $B^{2m}$ is orientable,
and the characteristic number $\int_B \omega_y$ is nonzero. In
\cite{FZ3} this characteristic number was called the  Weinstein
invariant associated to $y\in \ft$ and was computed explicitly in
terms of Chern and Pontrjagin numbers. For each adjoint orbit in
$\ft$ one obtains obstructions to fatness in terms of Chern and
Pontryagin numbers. The above discussion easily generalizes to fiber
bundles associated to principle bundles to obtain obstructions. We
call the metric on the fiber bundle a connection metric if the
fibers are totally geodesic. Some sample obstructions are:
\begin{thm*}
\begin{enumerate}
\item[(a)] (\cite{DR}) The only  $\Sph^3$ bundle over $\Sph^4$ that admits a fat connection metric is
the Hopf bundle $\Sph^3\to\Sph^7\to\Sph^4$.
  \item[(b)] \cite{FZ3} The only two  $\Sph^3$ sphere bundles  over $\CP^2$ that may
admit a fat connection metric are the complex sphere bundles
with $c_1^2=9$ and $c_2=1$ or $2$. In particular,
$T_1\CP^2\to\CP^2$ does not have a fat connection metric.
\item[(c)] \cite{FZ3} If a sphere bundle over $B^{2n}$
 admits a fat connection metric, then the Pontryagin numbers
 satisfy $\det(p_{j-i+1})_{1\leq i,j\leq n}\neq 0$.
\end{enumerate}
\end{thm*}
On the other hand, many of the examples in Section 2 are fat fiber
bundles, and the new example of positive curvature in Section 3 is a
fat bundle as well.

If one wants to achieve positive curvature on the total space, we
need to assume, in addition to the base having positive curvature,
that $G=\S^1$, $\SU(2)$ or $\SO(3)$. In \cite{CDR} a necessary and
sufficient condition for positive curvature of such metrics  was
given. The proof carries over immediately to the category of
orbifold principal bundles, which includes the case where the $G$
action on $P$ has only finite isotropy groups.

    \begin{thm}         \label{chavesderigas}
    (Chaves-Derdzinski-Rigas)
    A connection metric $g_t$ on an  orbifold $G$-principal bundle with $\dim G\le 3$
   has positive curvature, for $t$
    sufficiently small, if and only if
    $$
    \left(\nabla_x\Omega_u\right)(x,y)^2 < |i_x\Omega_u|^2 k_B(x,y),
    $$

    \noindent for all linearly independent horizontal vectors
     $x,y$ and $0\ne u\in\fg$.
    \end{thm}

    Here $k_B(x,y)=g(R_B(x,y)y,x)$ is the unnormalized sectional
    curvature and $i_x\Omega_u\ne 0 $
    is precisely the above fatness condition.

    All the known examples
    of principle bundles with positive curvature
    satisfy this condition (called hyperfatness
    in \cite{Zi2}). The new example in Section 3 also is such an
    $G=\S^3$ principle bundle, if one allows the action of $G$ to be
    almost free (the base in then an orbifold). All of the above
    easily carries over to obifold principle bundles as well.

    The condition is of course trivially satisfied if
    $\nabla\Omega=0$. This is equivalent to the metric on $B$ to
    be quaternionic K\"ahler. Thus a quaternionic K\"ahler manifold
    with \psc\ gives rise to a positively curved metric on the
    $\S^3$ or $\SO(3)$ principal bundle defined canonically by its
    structure. Unfortunately, if $\dim >4$, a positively curved
    quaternionic K\"ahler metric is isometric to $\HP^n$ only gives rise to $\Sph^{4n+3}(1)$
    on the total space of the principle bundle (Berger). And
    in dimension 4, the only smooth quaternionic K\"ahler metric
    (with positive scalar curvature) is isometric to $\Sph^4(1)$ and $\CP^2$
    giving rise to $\Sph^7(1)$ and the positively curved Aloff
    Wallach space $W_{1,1}$. Thus in the smooth category, it gives
    nothing new.

    But if the quotient is a quaternionic K\"ahler orbifold, there
    are many examples. In the case of the candidate $P_k$, the
    subgroup $\S^3\times\{e\}\subset \S^3\times\S^3$ acts almost
    freely and the quotient is an orbifold homeomorphic to $\Sph^4$.
    Similarly, for $Q_k$, the subgroup $\S^3\times\{e\}\subset \S^3\times\S^3$ acts almost
    freely with quotient  an orbifold homeomorphic to $\CP^2$. On
    these two orbifolds, Hitchin \cite{Hi} constructed a quaternionic K\"ahler metric
    with positive scalar curvature, and in \cite{GWZ} it was shown
    that the total space of the  canonically defined
    principal $\S^3$ resp. $\SO(3)$ principal bundle is
    equivariantly diffeomorphic to $P_k$ resp. $Q_k$. In \cite{Zi3}
    it was shown that the Hitchin metrics have a large open set on which the curvature is positive,
     but not quite everywhere. Nevertheless, in
     \cite{GWZ} and \cite{De} this Hitchin metric was the starting
     point. In  \cite{GWZ} a connection metric was constructed on $P_2$ by defining the
     metric piecewise
     via  low degree polynomials, both for the metric on the base, and the principle
     connection. The metric is indeed very close to the Hitchin
     metric. In \cite{De} the quaternionic K\"ahler Hitchin metric on
     $\Sph^4$ was deformed, but the the principle connection stayed
     the same. The metric on $\Sph^4$ was then approximated by
     polynomials in order to show the new metric has positive
     curvature.

     In both cases, the proof that the polynomial metric has positive
     curvature, was carried out by using a method due to Thorpe
     (and a small modification of it in \cite{De}).

     We finish by describing this Thorpe method since it is not so
     well known, but very powerful (see \cite{Th1}, \cite{Th2} and also \cite{Pu}).
     In fact, one of the problems of
     finding new examples is that after
     constructing a metric,  showing that it has positive curvature
     is difficult. Even the linear algebra problem for a
     curvature tensor on a vector space is  highly non-trivial! Here
     is where the Thorpe method helps.

Let  $V$ be a vector space with an inner product and $R$ a 3-1
tensor which satisfies the usual symmetry properties of a curvature
tensor. We can regard $R$ as a linear map
$$\hat{R}\colon\Lambda^2 V\to\Lambda^2 V ,$$ which, with respect to
the natural induced inner product on $\Lambda^2 V$, becomes a
symmetric endomorphism. The sectional curvature is then given by:
$$\sec(v,w)=\ml \hat{R}(v\wedge w),v\wedge w\mr
$$
if $v,w$ is an orthonormal basis of the 2-plane they span.

\smallskip

If $\hR$ is positive definite, the sectional curvature is clearly
positive as well. But this condition is exceedingly strong since a
manifold with $\hR>0$ is covered by a sphere \cite{BW}. But one can
modify the curvature operator by using a 4-form $\eta\in
\Lambda^4(V)$. It induces another symmetric endomorphism
$\hat{\eta}\colon\Lambda^2 V\to\Lambda^2 V$ via $\ml \hat{\eta}(x
\wedge y),z\wedge w\mr =\eta(x,y,z,w)$. We can then consider the
modified curvature operator $\hR_\eta=\hR+\hat{\eta}$. It satisfies
all symmetries of a curvature tensor, except for the Bianchi
identity. Clearly $\hR$ and $\hR_\eta$ have the same sectional
curvature since
$$\ml \hR_\eta(v\wedge w),v\wedge w\mr=
\ml \hR(v\wedge w),v\wedge w\mr + \eta(v,w,v,w)=\sec(v,w)$$

If we can thus find a 4-form $\eta$ with $\hR_\eta>0$, the sectional
curvature is positive.  Thorpe showed \cite{Th2} that in dimension
4,  one can always find a 4-form such that the smallest eigenvalue
of $\hR_\eta$ is also the minimum of the sectional curvature, and
similarly a possibly different 4-form such that the largest
eigenvalue of $\hat{R_\eta}$ is the maximum of the sectional
curvature. Indeed, if an eigenvector $\omega$ to the largest
eigenvalue of $\hat{R}$ is decomposable, the eigenvalue is clearly a
sectional curvature. If it is not decomposable, then
$\omega\wedge\omega\ne 0$ and one easily sees that $\hat{R_\eta}$
with $\eta=\omega\wedge\omega$ has a larger eigenvalue.

This is not the case anymore in dimension bigger than 4 \cite{Zo}.
Nevertheless this can be an efficient method to estimate the
sectional curvature of a metric. In fact, P\"uttmann \cite{Pu} used
this to compute the maximum and minimum of the sectional curvature
of all positively curved homogeneous spaces, which are not spheres.
It is peculiar to note though that this method does not work to
determine which homogeneous metrics on $\Sph^7$ have positive
curvature, see \cite{VZ1}.

\bigskip

\providecommand{\bysame}{\leavevmode\hbox
to3em{\hrulefill}\thinspace}

\end{document}